\RequirePackage[l2tabu, orthodox]{nag}

\documentclass[12pt]{amsart}
\usepackage{fullpage,url,amssymb,enumerate,colonequals}
\usepackage[all]{xy} 
\usepackage{mathrsfs} 
\usepackage[dvipsnames]{xcolor}
\usepackage{tikz}
\usetikzlibrary{shapes,snakes,arrows,positioning}
\usepackage{tikz-cd}
\usepackage{graphicx}
\usepackage{subcaption}
\usepackage{placeins}
\usepackage{todonotes}
\usepackage{tipa}

\usepackage[OT2,T1]{fontenc}
\DeclareSymbolFont{cyrletters}{OT2}{wncyr}{m}{n}
\DeclareMathSymbol{\Sha}{\mathalpha}{cyrletters}{"58}

\usepackage{color}

\newcommand{\defn}[1]{\emph{#1}} 

\newcommand{\F}{\mathbb{F}}

\newcommand{\PP}{\mathbb{P}}
\newcommand{\Q}{\mathbb{Q}}

\newcommand{\Z}{\mathbb{Z}}

\newcommand{\Zhat}{{\widehat{\Z}}}

\DeclareMathOperator{\Aut}{Aut}

\DeclareMathOperator{\End}{End}

\DeclareMathOperator{\Gal}{Gal}

\DeclareMathOperator{\img}{img}

\DeclareMathOperator{\Jac}{Jac}

\DeclareMathOperator{\ST}{ST}

\newcommand{\GL}{\operatorname{GL}}

\newcommand{\PSL}{\operatorname{PSL}}
\newcommand{\SL}{\operatorname{SL}}

\newcommand{\Kgal}{K^{\textrm{gal}}}

\newcommand{\wt}{\widetilde}

\newcommand{\gplabel}[1]{\href{https://www.lmfdb.org/Groups/Abstract/#1}{\texttt{#1}}}
\newcommand{\sublabel}[1]{\href{https://beta.lmfdb.org/Groups/Abstract/sub/#1}{\texttt{#1}}}

\newcommand{\GAP}{\textsf{GAP}}

\newcommand{\SageMath}{\textsf{SageMath}}

\newtheorem{theorem}{Theorem}[section]

\theoremstyle{definition}
\newtheorem{definition}[theorem]{Definition}

\newtheorem{example}[theorem]{Example}

\theoremstyle{remark}
\newtheorem{remark}[theorem]{Remark}

\makeatletter
\g@addto@macro\bfseries{\boldmath} 
\makeatother

\usepackage{microtype}

\usepackage[hyperfootnotes=false, colorlinks, linkcolor={blue}, citecolor={magenta}, filecolor={blue}, urlcolor={magenta}, plainpages=false, pdfpagelabels, pagebackref]{hyperref}

\begin{document}

\title{Creating a dynamic database of finite groups}

\date{\today}

\author[Combes]{Lewis Combes}
\address{Lewis Combes, School of Mathematics and Statistics, University of Sheffield, UK, S3 7RH}
\email{\url{lmcombes1@sheffield.ac.uk}}

\author[Jones]{John Jones}
\address{John Jones, Department of Mathematics and Statistical Sciences, P.O.\ Box 871804, Arizona State University, Tempe, AZ, 85287}
\email{\url{jj@asu.edu}}

\author[Paulhus]{Jennifer Paulhus}
\address{Jennifer Paulhus, Department of Mathematics and Statistics, Mount Holyoke College, Mount Holyoke,  MA, 01075}
\email{\url{jpaulhus@mtholyoke.edu}}

\author[Roe]{David Roe}
\address{
  David Roe,
  Massachusetts Institute of Technology,
  Department of Mathematics,
  77 Massachusetts Ave.,
  Bldg. 2-336
  Cambridge,
  MA 02139,
  United States of America
}
\email{\url{roed@mit.edu}}

\author[Roy]{Manami Roy}\address{Manami Roy, Department of Mathematics, Lafayette College, 730 High St, Easton, PA 18042}
\email{\url{royma@lafayette.edu}}

\author[Schiavone]{Sam Schiavone}
\address{
  Sam Schiavone,
  Massachusetts Institute of Technology,
  Department of Mathematics,
  77 Massachusetts Ave.,
  Bldg. 2-336
  Cambridge,
  MA 02139,
  United States of America
}
\email{\url{sam.schiavone@gmail.com}}


\begin{abstract}
A database of abstract groups has been added to the \emph{L-functions and Modular Forms Database} (LMFDB), available at \url{https://www.lmfdb.org/Groups/Abstract/}.  We discuss the functionality of the database and what makes it distinct from other available databases of abstract groups. We describe solutions to mathematical problems we encountered while creating the database, as well as connections between the abstract groups database and other collections of objects in the LMFDB.
\end{abstract}

\maketitle

\tableofcontents

\section{Introduction}\label{S:Intro}

Finite groups have played a profound role in mathematics for close to two centuries and, almost since their inception, mathematicians have asked classification questions about them. The quest to classify finite simple groups took up most of the $20^{\text{th}}$ century.
In the last 30 years as computational power became prominent, various databases of groups have been developed.  The \textit{ATLAS of Finite Groups} which contains information on interesting groups, notably character tables,  began as a book \cite{atlasbook} and is now available online \cite{atlasonline}.  The computer algebra programs \textsf{GAP} \cite{GAP, gap-survey} and \textsf{Magma} \cite{Magma} both have complete databases of groups of order up to 2000 (except for order $1024$, for which there are almost $49.5$ billion distinct isomorphism classes).  These computer algebra programs also include other databases such as almost simple and perfect groups (up to orders $1.6\cdot 10^7$ and $50{,}000$, respectively) and transitive groups up to degree $48$ \cite{gap-trans}, \cite{magma-trans}. The \textit{GroupNames} project \cite{groupnames} includes groups up to order 500 (skipping orders 256 and 384) as well as additional information such as  lattice of subgroups and character tables for the groups.

\subsection{Highlights of the database}
As we note above, there are several other databases of groups already available, so why do we need yet another one? Some of the useful features of our database are the following.

\begin{itemize}
    \item
        \textit{The database is dynamically searchable.} We have precomputed and stored many invariants and properties of the groups in the collection, and created a search interface which allows users to search on these stored values. For instance, we can quickly determine the $8$ groups $G$ of order $1440$ which satisfy the short exact sequence $$1 \to D_{10} \to G \to F_9 \to 1.$$ We can immediately determine that the smallest group to have a  character of degree $7$ is a group of order $56$. And in a matter of seconds we can find all 129,069 subgroups of groups of order between $1$ and $255$ which are non-normal, non-abelian, non-maximal, and non-trivial proper. The list of all such subgroups (a 23 megabyte file) can be downloaded in a few seconds. Each of these searches would be time-consuming and computationally intensive in computer algebra systems like \textsf{Magma} or \textsf{GAP}.


    \item
        \textit{The database is free and easy to use.} Use of our database requires no prior knowledge of any programming language or computer algebra package, making it accessible to a broad audience, from students just learning the basics of group theory, to experts using groups in their research. While most of our computations were done in \textsf{Magma}, which is closed-source and requires purchase of a license, the data we have computed is now publicly available on the internet through the \textit{L-functions and Modular Forms Database} (LMFDB).

    \item
        \textit{The database fosters connections between other collections of objects in the LMFDB.} The LMFDB is a huge international collaboration to  collect, curate, and connect computational mathematical work in number theory, particularly as it relates to the Langlands program.  Groups are attached to many mathematical objects throughout the LMFDB, from Sato-Tate groups to automorphism groups of curves,  and our database of groups facilitates more connections between these various collections. See Figure~\ref{connectiondiagram} in Subsection \ref{S:connections} for a diagram illustrating some of these connections.

    \item
        \textit{The database aggregates groups drawn from a variety of sources.} The database is flexible in that it can accommodate groups of various types, such as polycyclic, permutation, and matrix groups, drawn from multiple sources. See section \ref{S:data} for a list of sources used to compute the data.

    \item
        \textit{The database supports searches on subgroups.} For each group in the database, we have computed and stored information on its subgroups. Users can search for all groups containing a subgroup with a given property; for instance, one can easily find the 12 groups of order 96 that contain a normal subgroup isomorphic to $A_4$. (See Subsection \ref{SS:subgroup} for more details on the subgroup information stored.)

\end{itemize}

\subsection{Structure of the paper}
We begin in Section \ref{S:data} with information about sources for the groups currently included in the database, and how the data about them is generated. In Section \ref{S:features} we describe the web interface and demonstrate a number of features of the database. In the process of designing the database and computing the corresponding data, we encountered a series of challenges which required mathematical solutions. The next two sections discuss these challenges: Section \ref{S:presentations} describes the method by which we computed presentations for solvable groups and Section \ref{S:labels} details the procedures we used to deterministically label subgroups, conjugacy classes, and characters. In Section \ref{S:examples} we give examples 
of connections between the abstract groups database to other collections of objects in the LMFDB.

 \subsection{Acknowledgments} 
      We owe a large debt of gratitude to Tim Dokchitser, whose \textit{GroupNames} website  was the initial inspiration for our database and who generously shared code and database expertise with us (as well as giving us permission to re-post many of the definitions on his page as  knowls on the LMFDB). We are also appreciative of the advice and many suggestions Andrew Sutherland gave us.  We are grateful for various conversations and code from Michael Bush, Derek Holt, Alexander Hulpke, Bjorn Poonen, and  David Roberts. We thank the American Institute of Mathematics for hosting a mini-workshop where much of the initial planning happened, and the Institute for Computational and Experimental Research in Mathematics for a workshop where the project took off.

Roe and Schiavone were supported by the Simons Collaboration in Arithmetic Geometry, Number Theory, and Computation via Simons Foundation grant 550033. Paulhus was partially supported by a Frank and Roberta Furbush Scholarship from Grinnell College. Combes was supported by the Engineering and Physical Sciences Research Council grant EP/R513313/1. 

\section{Data overview}\label{S:data}

\subsection{Sources of data}
The database currently contains roughly 540,000 groups, 275 million subgroups and 40 million of their irreducible complex characters.


The code used to generate the data is written in \textsf{Magma} \cite{Magma} and \textsf{Python} \cite{python}, and it may be found at the GitHub repository \url{https://github.com/roed314/FiniteGroups}.  The database, as with the rest of the LMFDB, uses PostgreSQL as its database management system.  Computations were carried out via a combination of an AMD Epyc 7713 server with 256 2.0GHz cores and 2TB of RAM, and a distributed computation on Google Compute Engine, with a substantial use of GNU parallel \cite{parallel} in both cases.

\begin{table}[htb]
\begin{tabular}{r | c c c c c c}
 Source & Total & Solvable & Perm. & Matrix & OptimizedPC & MinPerm \\ \hline 
 Small Groups & 257936 & 257500 & 257746 & 68042 & 257500 & 257746 \\
 Transitive Groups & 235919 & 211279 & 235919 & 218 & 14499 & 161656 \\
 Intransitive Groups & 5444 & 2739 & 5444 & 16 & 2330 & 5378 \\
 Classical Lie Type & 2201 & 2 & 1509 & 2201 & 0 & 1509 \\
 CARAT & 189 & 185 & 186 & 189 & 174 & 186 \\
 $\text{GL}_n(\mathbb{F}_q)$ Subgroups & 3018 & 2456 & 3000 & 3018 & 2397 & 3000 \\
 $\text{GL}_2(\mathbb{Z}/N)$ Subgroups & 29771 & 28819 & 25319 & 29771 & 24323 & 25254 \\
 Perfect & 123 & 0 & 123 & 1 & 0 & 123 \\
Chevalley & 13 & 0 & 7 & 13 & 0 & 7 \\
 Sporadic & 9 & 0 & 9 & 7 & 0 & 9 \\
 Small Group Auto. & 283 & 283 & 283 & 0 & 73 & 111 \\
 Transitive Group Auto. & 498 & 438 & 498 & 0 & 51 & 201 \\
 Auto. Groups of Curves & 530 & 527 & 530 & 0 & 526 & 530 \\
\end{tabular}
\caption{Number of groups by source. See Section \ref{S:data} for more information. \label{davidsbigtable}}
\end{table}

We compute and store data on groups from various sources. Counts of groups currently in the database from these different sources are given in Table~\ref{davidsbigtable}. The ``Total'' column represents the total number of groups in the database from that source, excluding those already shown in a previous line. The  ``Solvable'', ``Perm.'', and ``Matrix'' columns give the number of groups from each source that are solvable, the number of groups for which we store a permutation representation, and the number of groups for which we store a matrix representation (over $\mathbb{Z}$, $\mathbb{F}_p$, $\mathbb{F}_q$, or $\mathbb{Z}/N$), respectively.  The ``OptimizedPC'' column  counts how many groups from that particular source include an optimized polycyclic presentation (see Section \ref{S:presentations}), meaning a presentation guaranteed to have a minimal number of generators, and ``MinPerm'' gives the number of groups for which we know a minimal degree permutation representation.

Roughly half of the data is initiated from the Small Groups library in \textsf{Magma} \cite{SmGp-1,SmGp-2,SmGp-3,SmGp-4,SmGp-5,SmGp-6,SmGp-7,SmGp-8,SmGp-9,SmGp-9addendum,SmGp-10}. All groups of order up to 2000, except those whose order is larger than 500 and divisible by 128, are included. In the rest of the paper, these are the groups we will call ``small groups''. Most of the remaining groups come from the transitive groups database \cite{TR-1,TR-2,TR-3}, from which we include all groups of degree up to 47, except those of degree 32, with orders between 512 and 40 billion.  Note that our notion of equivalence in this database is abstract isomorphism rather than conjugacy within $S_n$; the work to divide the transitive groups into isomorphism classes is described briefly in Subsection \ref{SS:isomorphism}.  In addition to the databases of small groups and transitive groups, we use the following sources:

\begin{itemize}
\item classical Lie groups up to particular bounds, and Chevalley groups that don't already show up as classical Lie groups \cite{Lie-1,Lie-2,Lie-3}; 
\item additional intransitive groups which are subgroups of $S_{15}$ not currently in the  database, computed using the \textsf{Magma} \texttt{Subgroup} command;
\item   all integer matrix groups of dimension up to 6 from CARAT \cite{Carat2.1b1};
\item subgroups of $\text{GL}_n(\mathbb{F}_q)$ computed via \textsf{Magma} for $n=2$ ($q < 1000$), $n=3$ ($q < 16$), $n=4$ ($q < 7$) and $n=5$ ($q=2$);
\item subgroups of $\GL_2(\Z/N\Z)$ for $N$ up to $125$ (skipping $80$, $96$, $104$, $112$ and $120$);
\item  all perfect groups of order up to 50,000 from the corresponding database in \textsf{Magma}; 
\item the sporadic groups $\mathrm{J}_1$, $\mathrm{J}_2$, $\mathrm{HS}$, $\mathrm{J}_3$, $\mathrm{McL}$, $\mathrm{He}$, $\mathrm{Ru}$, $\mathrm{Co}_3$, and $\mathrm{Co}_2$ with permutation and matrix representations taken from the Atlas of Sporadic Groups \cite{atlasonline};
\item automorphism  groups of curves  up to genus 48 \cite{breuer}, and ``large''  automorphism groups (of order $>4(g-1)$) of curves up to genus 101 \cite{conder}.
\end{itemize}

Since the automorphism groups of small groups can be much larger than the groups themselves, we also include all automorphism groups of groups of order up to $255$ and of transitive groups of degree up to 23.

\subsection{Basic attributes computed for every group}

For each group, we compute as many attributes as possible.  Some use commands directly from \textsf{Magma}, such as determining if a group is abelian, while others require special functions we wrote (for example, determining if a group is metacyclic).  We compute a ``reasonable'' presentation for the group (see Section \ref{S:presentations}), a minimal degree permutation representation, the lattice of subgroups up to automorphisms of the group, and up to conjugacy when possible (Section \ref{SS:subgroup}), and conjugacy classes and the character tables when feasible. We also determine special subgroups such as the center, commutator, Frattini, and Fitting subgroups, as well as various series for the group.

\subsection{Challenging attributes computed for many groups}

Many of our quantities of interest become more difficult to compute as the size or degree of the group increases. Our general approach is to try to compute everything for all groups, setting appropriate timeouts and omitting data if the computation does not finish or if we encounter errors in \textsf{Magma}. We have developed additional code for some problems where \textsf{Magma}'s built-in methods were insufficient for our purposes; notable examples include computing quotients of permutation groups and isomorphism testing (discussed in Subsection \ref{SS:isomorphism}).

\subsection{Computing subgroup lattices}

The subgroup lattice is an example of a particularly challenging attribute to compute. For each group, we attempt to compute the full lattice of subgroups, but there are many groups in our database where it is infeasible to compute subgroups up to conjugacy.  For example, $C_2^{10}$ has 229,755,605 subgroups, none of which are conjugate. However, it only has 11 classes of subgroups up to automorphism.  We implemented the computation of the subgroup lattice up to automorphism and improved on \textsf{Magma}'s built-in subgroup lattice methods up to conjugacy.  In some cases, such as $S_{47}$, working up to automorphism does not sufficiently reduce the quantity of subgroups, so we restrict our attention to normal subgroups and their complements, maximal subgroups, Sylow subgroups, subgroups with small index and/or trivial core due to their importance in permutation representations.

\textsf{Magma} has built-in methods for computing both the subgroup lattice and the list of subgroups up to conjugacy (without inclusions).  We found that the second was much faster than the first, and that we could recover the inclusions more quickly by computing a vector counting the intersections with each conjugacy class: given a group $G$ and subgroups $H_1 \subseteq H_2$, then the vector of counts for $H_2$ dominates that for $H_1$, dramatically decreasing the number of calls to \texttt{IsConjugateSubgroup}.  When there are many subgroups up to conjugacy (e.g.,\ abelian $p$-groups), we instead compute subgroups up to automorphism.  The first algorithm for doing so uses the holomorph of $G$, i.e., the semidirect product $G \rtimes \Aut(G)$. Conjugacy within the holomorph translates to automorphism in the group itself, and a code snippet kindly provided by Derek Holt using lifting through an elementary abelian series allowed for the computation of representatives of subgroups up to automorphism in the solvable case without computing the list up to conjugacy.

However, the holomorph is implemented in \textsf{Magma} as a permutation group of degree equal to the order of $G$, and as the size of $G$ grows, computations using the holomorph bog down.  For non-solvable $G$ and $G$ with order at least $5000$, we switch to a graph theoretic algorithm.  First we compute a list of subgroups up to conjugacy, together with a list of automorphisms that generate the outer automorphism group of $G$.  Taking subgroups up to conjugacy as vertices, we add an edge between $H_1$ and $H_2$ if there is an outer generator $\sigma$ with $\sigma(H_1)$ conjugate to $H_2$.  The components of the resulting graph give the subgroups up to automorphism. 


\subsection{The isomorphism problem} \label{SS:isomorphism}

When adding transitive permutation groups to our abstract groups database, we had to ensure that groups with multiple transitive permutation representations were only added once each.  For small groups, \textsf{Magma} will compute its small group identification number.  In a larger range, the LMFDB's transitive group section already had the needed information.  However, there were multiple cases where neither approach sufficed.

Running isomorphism tests can be time consuming, so we used the following strategy.  We took multiple isomorphism invariants of a group which are quick to compute and combined them into a hash.  This could be precomputed for large numbers of groups.  Groups with different hashes were clearly non-isomorphic.  If groups produced identical hashes, we then had \textsf{Magma} perform a slower, but conclusive test for isomorphism.  Our hash was highly effective in distinguishing non-isomorphic groups.  Full details will be given in the forthcoming article \cite{roehash}.



\section{Features of the database}\label{S:features}

The database is viewable through the LMFDB website at \url{https://www.lmfdb.org/Groups/Abstract/}. In this section we describe how to navigate the web interface, and highlight some particular features.

\subsection{Features of the search interface}
Groups can be searched through numerous query types, for example: the isomorphism type of a group's automorphism group, commutator, center, abelianization, Frattini subgroup, etc.; the order of any of these and the group itself; and many boolean properties such as nilpotency, simplicity, and solubility. With the search interface, it is easy to quickly answer questions such as \textit{``What are the nilpotent groups of order 36?''} and \textit{``Are there any groups of order 256 with abelianization $C_2 \times C_2$?''} In this way, the database provides a useful service to researchers looking for groups with particular properties, as well as students coming to grips with the basics of group theory. It is also possible to search for groups whose orders factor in a particular way. For example, groups whose order is of the form $p^3q^2r$ for primes $p,q,r$ can be queried with the string $[3,2,1]$. 

There is a \href{https://www.lmfdb.org/Groups/Abstract/interesting}{curated list} of some interesting groups that one can quickly access, and one can view the home page of a group \href{https://www.lmfdb.org/Groups/Abstract/random}{picked randomly} from the database as well. Each abstract group has a unique label attached to it (see Section~\ref{S:labels} for more details), and one may search for a specific abstract group using its label or by a name such as ``S5''.

\subsection{Features of a group page} \label{SS:group_page}
On individual group pages, information is organized roughly by topic. Constructions of the group via a presentation (see Section \ref{S:presentations}), and as direct, semidirect, and non-split products are given when possible, as well as representations as a matrix or permutation group. Homological information like the Schur multiplier and commutator length are provided. 
Special subgroups like the center, commutator and socle are displayed.  For many groups we provide several subgroup diagrams; one can toggle whether subgroups are shown up to conjugacy or up to automorphism (see Subsection~\ref{SS:subgroup}), and independently whether to show all such classes of subgroups or just normal subgroups. 
We also display the derived series, chief series,  and upper and lower central series of the group. Groups for which the given group is a maximal subgroup or a maximal quotient are listed under the Supergroups section.

\begin{figure}[ht]
    \centering
    \includegraphics[width=0.6\textwidth]{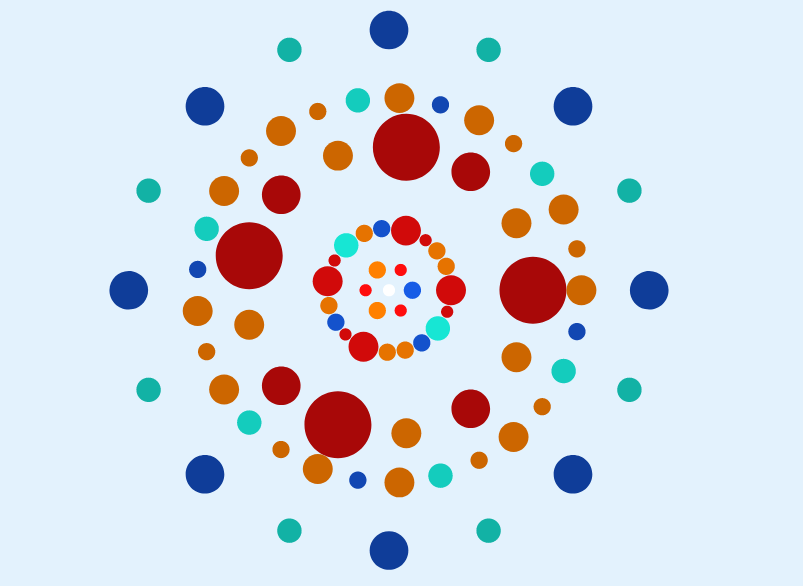}
    \caption{Group picture of \href{https://www.lmfdb.org/Groups/Abstract/480.60}{\texttt{480.60}} \label{grouppic}}
\end{figure}

Groups are also represented pictorially, in line with other sections of the LMFDB. Each disc in a group picture represents a conjugacy class of elements, arranged in an annulus around the middle class (of the identity), with distance according to the number of prime factors in the class's order. A disc's size and color are chosen based on its order. The decisions on how to display a group in picture form are essentially a matter of taste, with some group pictures showing striking symmetries within the conjugacy classes, such as \gplabel{480.60} in Figure~\ref{grouppic}.

For many groups we provide character tables for both complex and rational characters.  When the data has been computed but the table is too large, the user is provided with a link to display the table.  Rows correspond to irreducible characters, columns correspond to conjugacy classes for the complex character table and divisions for rational character table. (For more details on conjugacy classes and divisions, see Subsection \ref{SS:CC_label}.)  For ease of terminology, in this discussion we refer to columns as corresponding to ``classes'' to cover both cases simultaneously.

The headers for classes expand to show basic information such as the class's size and centralizer, as well as a representative for the class. As is standard, the row(s) above the table contain power map information on classes.  For each prime $p$ dividing the order of the group and class represented by an element $g$, the entry above that column in the row labeled ``pP'' gives the label for the class of $g^p$.  For each complex character, in the 2nd column of the table we designate whether the representation is real, complex, or quaternionic.

\begin{remark}
   For some groups which lie outside the range of the database, we still generate webpages that display some basic information about the group by using on-the-fly calculations in \GAP{} and \SageMath. Two types of groups that are handled this way are groups in \GAP’s Small Groups database that are outside the LMFDB collection, such as \gplabel{512.402873}, and large abelian groups, such as \href{https://www.lmfdb.org/Groups/Abstract/ab/12_3.6}{$C_{12}^3\times C_6$}.
\end{remark}

\begin{remark}
    Statistics of the abstract groups database can be found at \url{https://www.lmfdb.org/Groups/Abstract/stats}. We have computed many statistics based on the order factorization type of a group, i.e., whether the order is a prime, a product of two distinct primes, a power of one prime times another prime, etc. For instance, one can see the distribution of the nilpotency class of groups with order $p^k$ where $p$ is a prime and $k \in \{3, \ldots, 6\}$.
\end{remark}

\subsection{Subgroups}\label{SS:subgroup}

One of the main features of our work is a companion database of subgroups.  Subgroups have their own  search functionality and labeling system (see Subsection \ref{SS:subgroup_label}). Given a group $H$, it is possible to search the database for all groups $G$ containing $H$ as a subgroup, subject to numerical constraints such as the index of $H$ in $G$ and the order of $G$. Cyclic subgroups, normal subgroups, and Sylow and Hall subgroups of a group can all be searched for. Each subgroup also has its own page. These pages give information about the ambient group and the quotient group structure (when normal). Related subgroups such as the centralizer, normalizer, normal closure and core of the subgroup are given.    

\begin{figure}[htb]
    \centering
    \includegraphics[width=1\textwidth]{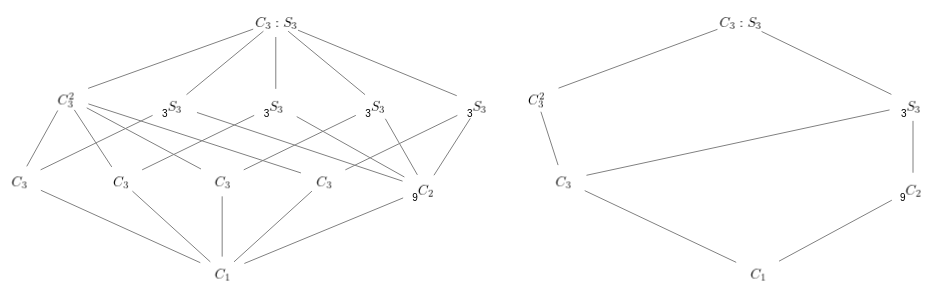}
    \caption{Subgroup diagrams of \gplabel{18.4} up to conjugacy and autjugacy, as they appear on the LMFDB}
    \label{fig:subgroup_diagrams}
\end{figure}

The subgroups database allows us to display the subgroup lattice for many groups.  This lattice is algorithmically very useful, allowing the computation of the rank of a group, its expressions as a semidirect product, etc. The subgroup lattice can be viewed on a group page when it has fewer than 100 subgroup classes; see Figure~\ref{fig:subgroup_diagrams} for an example. Inclusion is indicated with lines, and clicking a particular subgroup reveals information about it, e.g., whether it is maximal or solvable, its normalizer, its core, etc. Subgroups are separated into levels vertically based on their orders or by the number of prime factors in their order (the user can toggle between these two options), and can be rearranged by clicking and dragging. When the diagram is too large to conveniently display on the page, the full diagram can often be viewed on a separate linked page.  

In addition, we provide profiles of the subgroups up to automorphism and up to conjugation.  For each order of subgroup, the profile lists the isomorphism types and multiplicities of classes of subgroups with that type.  For example, for the elementary abelian group $C_2^5$, the profile for subgroups up to conjugation would list $155$ conjugacy classes of subgroups of order $8$ (all isomorphic to $C_2^3$), and one class of subgroups of the same order up to automorphism.

\section{Computing presentations for solvable groups}\label{S:presentations}

As mentioned in Subsection \ref{SS:group_page}, each group page contains various constructions of the group. In the case of solvable groups we give particularly nice presentations, exploiting the fact that solvable groups are polycyclic.
A group $G$ is \textit{polycyclic} if it has a descending chain of subgroups
\begin{align*}
    G = G_1 \geq G_2 \geq \cdots \geq G_{n+1} = 1
\end{align*}
such that each $G_{i+1}$ is normal in $G_i$, and the quotient $G_i/G_{i+1}$ is cyclic. A polycyclic presentation is a special type of presentation for the group $G$ that takes advantage of this descending chain of subgroups; for details, see \cite[Definition 8.7]{Handbook}.

However, the polycyclic presentations in \textsf{Magma} use one generator for each prime dividing the order of the group (counted with multiplicity).  So, for example, the cyclic group of order $2^2 \cdot 3 \cdot 5$ would have a presentation with $4$ generators instead of the simpler one using a single generator.
Here we describe a process for finding a more human-readable polycyclic presentation for a solvable group by using chains of subgroups.  One can pass between polycyclic presentations and subnormal filtrations with cyclic quotients:
\begin{enumerate}
\item given a polycyclic presentation with generators $\{g_1, \dots, g_m\}$, 
\[
1 \le H_1 \le \dots \le H_m = G
\]
will be a filtration of $G$ with cyclic quotients, where $H_i = \langle g_1, \dots, g_i \rangle$;
\item given a filtration of the form above, any choice of elements $g_i$ generating $H_i / H_{i-1}$ will give generators for a polycyclic presentation.
\end{enumerate}
Define the \defn{relative order} of $g_i$ to be the order of the quotient $H_i / H_{i-1}$.  We approach the search for a presentation by first finding a filtration and then arbitrarily choosing generators.

 We first construct all minimal length chains of subgroups, each normal in the next, with cyclic relative quotients. We do this via a ``top-down'' approach, building them down from the top in layers until reaching the trivial subgroup.  This process guarantees a presentation with a minimal number of generators; we then compare all such minimal presentations according to the following criteria.
 \begin{enumerate}
    \item Maximize the number of generators with order equal to their relative order.
    \item Maximize the number of generators that commute with each other.
    \item Aim for relative orders that are non-increasing.
    \item Conjugacy relations should be ``deeper'' (within smaller groups in the filtration).
\end{enumerate}
This process cannot feasibly be made canonical: at some points we make arbitrary choices, both for the filtration and for the choice of generators.

There are some groups for which the process of constructing and comparing all minimal length chains is too time-consuming.  In these cases we instead use two other algorithms, both of which run much faster, but may not give as good a presentation.  The first proceeds by greedily building the chains by starting with $H = G$ and iteratively choosing a random normal subgroup $K < H$ with $H/K$ cyclic.  The second algorithm simply refines the derived series by taking an abelian basis at each step in order to fill in the series to one with cyclic quotients.  

It is important to pick the presentation at the beginning of our computations for each solvable group, as this presentation impacts certain attributes of the group, such as how we represent elements of the group.


\section{Labels}\label{S:labels}

A foundational principle of the LMFDB is to display every object with a unique label, an identifier that carries mathematically relevant information which can ideally be computed in a deterministic way.  The primary role of the label is to serve as a permanent identifier for the object (so that it can be referenced in the mathematical literature for example).  Labels generally encode a sequence of invariants of an object, and provide a natural ordering.  Conversely, a sequential ordering of all objects of a type trivially allows one to label them, if by nothing else, by their position in the sequence.  We will talk of having labels and an ordering interchangeably.  When sorting ordered tuples, we always mean lexicographically, and treat complex numbers as ordered pairs of real numbers for this purpose.

We label small  groups as \texttt{N.i} corresponding to the \textsf{GAP} ID encoded as a string, where \texttt{N} is the order of the group and \texttt{i} distinguishes groups of the same order (as determined in \textsf{GAP}). If a group is not in the \textsf{GAP} Small Groups database, we replace \texttt{i} with an incrementing
letter code, assigning labels to groups as they are added to our database.


The existence of automorphisms precludes having a definitive labeling \emph{de novo}.  Consequently, the labels for conjugacy classes, characters, and subgroups depend on fixing an ordered list of generators, which we do, and a specific realization of the group, i.e., as a permutation group, a polycyclic group, or as a matrix group.  The latter is important since in each case, one can easily construct an injective set map from the group to $\Z$, which in turn orders elements of the group.
We then need to have a reproducible method for generating elements of a group.  In the applications below, we found that picking pseudo-random group elements worked efficiently; for details on how we generate such elements, see \url{https://www.lmfdb.org/knowledge/show/group.pseudo_random_elements}.

%
%

\subsection{Conjugacy classes and divisions} \label{SS:CC_label}


If $G$ is a group and $g_1, g_2\in G$, we say that $g_1$ and $g_2$ are in the same {\em division} if there exists $h\in G$ such that $h^{-1}g_1h=g_2^i$ for some $i\in\Z$ such that $\gcd(i,|g_1|)=1$.  Divisions are unions of conjugacy classes, and rational-valued characters are constant on divisions. In fact, the irreducible rational-valued characters form a basis for the space of functions that are constant on divisions, as complex characters do for class functions.  We say that a division is {\em maximal} if a representative generates a maximal cyclic subgroup of $G$.

Full details of how we label conjugacy classes and divisions are given at \url{https://www.lmfdb.org/knowledge/show/group.division_computing_labels}; here we describe some properties of the labeling.

Divisions are labeled with \texttt{nA} where $n$ is a positive integer giving the order of a representative element and \texttt{A} is a string made up of capital letters.  They are ordered first by the size of a conjugacy class within the division, and then by the number of conjugacy classes within a division.

Finding representatives for each division can be computationally expensive, so we focus first on maximal divisions.  From there, we can compute labels and representatives for other classes as follows.  Working through the maximal divisions in order, we use its chosen representative $g$ and compute the classes of elements $g^i$ with $0<i<|g|$ and $\gcd(i, |g|)>1$.  The order of the non-maximal divisions is simply the order in which they appear from this sequence (looping over divisions, and within that, powers of an element).

The label for a conjugacy class starts with the label of the division containing it.  If there is only one class within the division, the two labels are the same.  When there is more than one class in a division \texttt{nA}, the conjugacy classes are labeled \texttt{nAj} where $j$ is an integer.  The first class encountered will then be \texttt{nA1}.  Let $g$ be an element representing this class.  We then consider the conjugacy classes $[g^{-1}]$, $[g^2]$, $[g^{-2}]$, \ldots. These correspond to conjugacy class labels \texttt{nA-1}, \texttt{nA2}, \texttt{nA-2}, \ldots.  The label for the class is the first time it appears in this sequence.
Conjugacy classes which land in the same division have the same prefix, and will be grouped together in the complex character table for the group.

\subsection{Characters}

When feasible, we give information on irreducible rational-valued and complex-valued characters.  The group $\Gal(\overline{\Q}/\Q)$ acts on the complex characters.  The Galois orbits correspond to the rational characters, with each rational character being simply the sum of the complex characters in an orbit.  We note that a rational-valued character may not arise from a rational-valued representation.  The Schur index gives the smallest multiplier for a rational character so that the resulting character arises from a rational-valued representation.

Permutation representations of a group $G$ can be realized using permutation matrices, and are thus sums of irreducible representations.  Transitive permutation representations are classified by their degree $n$ and {\em $T$-number} $t$ \cite{conway-et-al}.   For a given irreducible representation, we can consider the ``first'' transitive permutation containing it to be the smallest pair $(n,t)$.

The labels for rational characters are of the form \texttt{G.na} where \texttt{G} is the label for the group, \texttt{n} is the degree, and \texttt{a} is a lower-case letter which acts as a counter.  The label for a complex character takes the form \texttt{G.nak} where \texttt{G.na} is the label for the corresponding rational character and \texttt{k} is positive integer which serves as a counter.  In the character tables, characters are sorted by their labels viewed as a tuple $(n,a)$ or $(n,a,k)$, respectively.  To completely determine the labels, we just need to know the ordering among characters with the same degree.

For rational characters, we sort by $(d,m,n,t)$ where $d$ is the degree (which is explicitly given in the label), $m$ is the size of the Galois orbit of complex characters giving rise to it, and $n$ and $t$ refer to the smallest containing permutation representation.  Any remaining ties are broken by sorting on the vector of character values using our ordering of conjugacy classes. Complex irreducible characters are given in the same order as their corresponding rational characters, with characters within a Galois orbit ordered by their vectors of values.

\subsection{Subgroups}\label{SS:subgroup_label}

Labels for subgroups take one of four forms, depending on whether they are computed up to conjugacy or up to automorphism, and whether or not the deterministic process described below succeeds. (If it fails or we do not have all subgroups of a given order, we use a fallback nondeterministic label.) The deterministic labels take the form $\texttt{n.i.m.a.c}$ (for subgroups up to conjugacy) or $\texttt{n.i.m.a}$ (for subgroups up to automorphism).  The piece $\texttt{n.i}$ is the label of the ambient group $G$, $\texttt{m}$ is the index of $H$ in $G$, and $\texttt{a}$ and $\texttt{c}$ are alpha-numeric identifiers that distinguish $H$ up to automorphism and conjugacy, respectively. To determine $\texttt{a}$ and $\texttt{c}$ we use the idea of Gassmann equivalence classes to order subgroups of the same index. For more on Gassmann equivalence and related notions see \cite{Sutherland-Equiv}, and for full details on the subgroup labeling scheme see \url{https://beta.lmfdb.org/knowledge/show/group.computing_subgroup_labels}.

\begin{definition}
Two subgroups $H_1$ and $H_2$ of $G$ are \defn{Gassmann equivalent} if they intersect each $G$-conjugacy class with the same cardinality. Equivalently, $H_1$ and $H_2$ have the same index $m$ and the permutation representations $\pi_{H_1}: G \to S_m$ and $\pi_{H_2}:G \to S_m$ have the same character (which counts the number of cosets fixed by each conjugacy class).
\end{definition}

In order to label Gassmann classes we first fix an ordering of the conjugacy classes of $G$ as in Section~\ref{SS:CC_label}. To each subgroup $H$ we then associate a vector of positive integers whose entries are the sizes of the intersections of $H$ with each conjugacy class. We then sort the subgroups using these vectors.  An analogous process works for Gassmann vectors up to automorphism, where we collect together $G$-conjugacy classes that are related by an automorphism of $G$.

If there are multiple subgroups up to conjugacy (resp., automorphism) within a given Gassmann class, we use the subgroup lattice to further order those remaining subgroups. We sort subgroups $H$ in the same Gassmann class by considering the collection of supergroups of $H$; if this does not suffice to distinguish Gassmann equivalent subgroups, we resort to computing permutation representation and their associated characters.

Once we have an ordering inside the Gassmann classes, we assign a letter label to the Gassmann class and concatenate a number from the ordering of subgroups inside the Gassmann class to create an alpha-numeric value which we assign to $\texttt{a}$ and (whenever we can compute up to conjugation) we similarly assign an alpha-numeric value to $\texttt{c}$.  For example, the label of $A_5$ inside $S_5$ (up to conjugacy) is $\sublabel{120.34.2.a1.a1}$. As another example, there are two classes of subgroups (up to automorphism) of $\operatorname{SL}(3,4)$ isomorphic to $\operatorname{GL}(2,4)$ with labels $\sublabel{60480.a.336.a1}$ and $\sublabel{60480.a.336.a2}$---these are Gassmann-equivalent but not related by an automorphism---and four classes of subgroups (up to automorphism) of $\operatorname{SL}(3,4)$ of order 8, each contained in a different Gassmann class: $\sublabel{60480.a.7560.a1}$ ($Q_8$), $\sublabel{60480.a.7560.b1}$ ($C_2 \times C_4$), $\sublabel{60480.a.7560.c1}$ ($D_4$), and $\sublabel{60480.a.7560.d1}$ ($C_2^3$).

In some cases, we store only a few subgroups of a given order (if we store the center of the group for example but not other subgroups of that order) or the code for computing the labels did not finish.  In these cases, we append an arbitrary letter code to the $\texttt{N.i.m}$ part of the label, with capital letters indicating that subgroups are stored up to conjugacy and lower case indicating subgroups up to automorphism.


\section{Connections and examples}\label{S:examples}


We anticipate that the group home pages will be a valuable tool for everyone from undergraduate students first learning abstract algebra to experts in the field. The layout of each page, the search features, and the knowls (the myriad of hyperlinks giving definitions of many of the words on each page) allow for exploration of many group theoretic concepts. Students can easily find examples of groups with common features like simple or solvable groups; they can explore Sylow subgroups (and more generally subgroup lattices); and they can learn about advanced topics such as series or character theory. Students can also explore relationships between different attributes of groups.  A quick search demonstrates examples of non-trivial groups that are perfect but not simple or $A$-groups (a group with all its Sylow subgroups abelian) that are not solvable.

For researchers, the abstract groups database is intimately connected with many other mathematical objects in the LMFDB.  We highlight some of those connections below.

\subsection{Connections with other objects} \label{S:connections}
\begin{figure}[htb]	
	\centering
	\scalebox{0.85}{
		\begin{tikzpicture}[node distance=1cm, auto]  
			\tikzset{
				mynode/.style={rectangle,rounded corners,draw=black, dashed, top color=white, bottom color=yellow!60,very thick, inner sep=0.3em, minimum size=5em, text centered},
				mynode1/.style={rectangle,rounded corners,draw=black, top color=white, bottom color=green!40,very thick, inner sep=0.2em, minimum size=3.5em, text centered},
				mynode2/.style={rectangle,draw=black, rounded corners, top color=white, bottom color=green!30,very thick, inner sep=0.3em, minimum size=3em, text centered},
				mynode21/.style={rectangle,draw=black ,rounded corners, dotted, top color=white, bottom color=black!15,very thick, inner sep=0.2em, minimum size=3em, text centered},
				mynode22/.style={circle,draw=black, top color=white, bottom color=blue!40,very thick, inner sep=0.2em, minimum size=3em, text centered},
				mynode3/.style={rectangle,rounded corners,draw=black, top color=white, bottom color=yellow!50,very thick, inner sep=0.3em, minimum size=2em, text centered},
				myarrow/.style={->, >=latex', shorten >=1pt, thick},
				myarrow1/.style={->, >=latex', shorten >=1pt, t	zigzag, hick},
				mylabel/.style={text width=7em, text centered} 
			}  
			\node[mynode1] (gal) {\ \ \ \ \ \ \begin{tabular}{c} \textbf{\href{https://www.lmfdb.org/GaloisGroup/}{Galois groups}} \\ (Transitive groups)\end{tabular}\ \ \ \ \ \ };
			\node[right=0.1cm of gal] (dummy) {}; 
			\node[left=0.1cm of gal] (dummy1) {};
			\node[right=4cm of gal] (dummy2) {};  
			\node[left=2cm of gal] (dummy3) {};  
			\node[mynode, below= 8cm of gal] (abs) {\ \ \ \ \ \ \textbf{\href{https://www.lmfdb.org/Groups/Abstract/}{Abstract groups}}  \ \ \ \ \ \ };
				\node[mynode2, below=4cm of dummy]  (NF) { \textbf{\href{https://www.lmfdb.org/NumberField/}{Number fields}}};
			\node[mynode2, below=4cm of dummy1] (padic) { \textbf{\href{https://www.lmfdb.org/padicField/}{$p$-adic fields}}};
			\node[mynode2, left=1.5cm of padic] (genus2) {\begin{tabular}{c} \textbf{\href{https://www.lmfdb.org/Genus2Curve/Q/}{Genus $2$ curves}} \\ \textbf{\href{https://www.lmfdb.org/Genus2Curve/Q/}{over $\Q$}} \end{tabular}};
			\node[mynode21, below= 6cm of dummy3] (Belyi) {\begin{tabular}{c}  \textbf{\href{https://www.lmfdb.org/Belyi/}{Belyi}} \\ \textbf{\href{https://www.lmfdb.org/Belyi/}{maps}}\end{tabular}};	
			\node[mynode2, below=1cm of dummy2] 	(mod) 	{\begin{tabular}{c} \textbf{\href{https://www.lmfdb.org/ModularForm/GL2/Q/holomorphic/}{Classical modular}}\\ \textbf{\href{https://www.lmfdb.org/ModularForm/GL2/Q/holomorphic/}{ forms of weight $1$}} \end{tabular}};
			\node[mynode2, below=4cm of mod] (Artin) {  \textbf{\href{https://www.lmfdb.org/ArtinRepresentation/}{Artin reps}}};
			\node[mynode2, right=3cm of abs] (HGC) {\begin{tabular}{c} \textbf{\href{https://www.lmfdb.org/HigherGenus/C/Aut/}{Higher genus}} \\ \textbf{\href{https://www.lmfdb.org/HigherGenus/C/Aut/}{curves}}\end{tabular}};
			\node[mynode2, left=3cm of abs] (EC) {\begin{tabular}{c} \textbf{\href{https://www.lmfdb.org/EllipticCurve/Q/}{Elliptic curves}}\\  \textbf{\href{https://www.lmfdb.org/EllipticCurve/Q/}{over $\Q$}} or  \textbf{\href{https://www.lmfdb.org/EllipticCurve/}{$\Q(\alpha)$ }} \end{tabular}};
			\node[below=5cm of abs] (dummy4) {}; 
			\node[mynode21, below= 0.1cm of dummy4] (Hyp) {\begin{tabular}{c} \textbf{\href{https://www.lmfdb.org/Motive/Hypergeometric/Q/}{Hypergeometric}} \\ \textbf{\href{https://www.lmfdb.org/Motive/Hypergeometric/Q/}{motives over $\Q$}}\end{tabular}};
			\node[mynode1, left= 2cm of  Hyp] (STgroups) {\begin{tabular}{c} \textbf{\href{https://www.lmfdb.org/SatoTateGroup/}{Sato-Tate groups}}\end{tabular}};
			\node[mynode21, right= 1cm of Hyp] (Lattices) {\begin{tabular}{c} \textbf{\href{https://www.lmfdb.org/Lattice/}{Lattices}}\end{tabular}};
			\node[mynode21, right= 2cm of Lattices] (ModCurves) {\begin{tabular}{c} \textbf{\href{https://www.lmfdb.org/ModularCurve/Q/}{Modular}} \\\textbf{\href{https://www.lmfdb.org/ModularCurve/Q/}{curves}}\end{tabular}};
			\draw[->, >=latex', shorten >=3pt, bend right=0,ultra thick] (gal.south) to (abs.north); 
			\draw[->, >=latex', shorten >=3pt, bend left=1,very thick] (STgroups.north) to node [left, swap] {\begin{tabular}{c} component \\groups\end{tabular}} (abs.south);
			\node[below=3cm of abs] (dummy5) {\begin{tabular}{c} monodromy \\ mod $\ell$\end{tabular}}; 
			\draw[-, >=latex', shorten >=3pt, bend right=0,very thick] (Hyp.north) to (dummy5.south);
			\draw[->, >=latex', shorten >=3pt, bend right=0,very thick] (dummy5.north) to (abs.south);    
			\draw[->, >=latex', shorten >=3pt, bend right=0,very thick] (Lattices.north) to node [right, swap] {\begin{tabular}{c} auto. \\ \hspace{0.01 in} groups\end{tabular}} (abs.south);  
			\draw[->, >=latex', shorten >=3pt, bend right=0,very thick] (ModCurves.north) to node [right, swap] {\begin{tabular}{c} \hspace{0.2 in}  subgroups of \\ \hspace{0.2 in} $\GL_2(\Z/N\Z)$\end{tabular}} (abs.south);  
			\draw[->, >=latex', shorten >=3pt, bend right=0,very thick] (HGC.west) to node [below, swap] {\begin{tabular}{c} auto.\\ groups \end{tabular}} (abs.east);  
			\draw[->, >=latex', shorten >=3pt, bend right=0,very thick] (EC.east) to node [below, swap] {\begin{tabular}{c} torsion \hspace{0.3in} \\ groups \hspace{0.3in} \end{tabular}} (abs.west);  
			\draw[->, >=latex', shorten >=3pt, bend right=5,very thick] (Belyi.south) to node [above, swap] {\begin{tabular}{c} auto.\\ groups \end{tabular}} (abs.north);  
			\draw[->, >=latex', shorten >=3pt, bend left=25,very thick] (Belyi.north) to node [left, swap] {\begin{tabular}{c} \ \ monodromy\\ groups \end{tabular}} (gal.west);  
			\draw[->, >=latex', shorten >=3pt, bend left=25,very thick] (padic.north) to node [left, swap] {\begin{tabular}{c} Galois\\ groups\ \ \end{tabular}} (gal.south);  
			\draw[->, >=latex', shorten >=3pt, bend left=15,very thick] (padic.south) to node [above, swap] {\begin{tabular}{c} \ \ \ \ inertia\\ \ \ \ \ \ \ groups\ \ \end{tabular}} (abs.north); 
			\draw[->, >=latex', shorten >=3pt, bend right=10,very thick] (Artin.west) to node[above, swap] {\begin{tabular}{c}  \hspace{0.1in}  projective image\vspace{0.1in} \end{tabular}} (abs.north);  
			\draw[->, >=latex', shorten >=3pt, bend left=5,very thick] (Artin.west) to node[below, swap] {\begin{tabular}{c}  \hspace{0.35in}  image \\ \end{tabular}} (abs.north);  
			\node[left=5.15cm of mod] (dummy6) {}; 
			\draw[<->, >=latex', shorten >=3pt, bend right=0, very thick] (mod.south) to (Artin.north);  
			\draw[->, >=latex', shorten >=3pt, bend right=12,very thick] (NF.north) to node [right, swap] {\begin{tabular}{c} Galois\\ \ \ groups \end{tabular}} (gal.south);  
			\draw[<-, >=latex', shorten >=3pt, bend left=10, thick] (NF.east) to (Artin.north);  
			\draw[->, >=latex', shorten >=3pt, bend left=10, thick] (NF.south) to  (Artin.west);   
			\node[above=2cm of genus2] (dummy8)  {\begin{tabular}{c} monodromy\\ groups \end{tabular}};
			\node[above=3.76cm of genus2] (dummy9){};
			\draw[-, >=latex', shorten >=3pt, bend left=0,very thick] (genus2.north) to (dummy8.south);   
			\draw[->, >=latex', shorten >=3pt, bend left=0,very thick] (dummy8.north) to (dummy9.north) to (gal.west);  
			\node[below=1.1cm of genus2] (dummy11)  {\begin{tabular}{c} auto.\\ groups \end{tabular}};
			\node[below=2cm of genus2] (dummy12){};
			\draw[-, >=latex', shorten >=3pt, bend left=0,very thick] (genus2.south) to (dummy11.north);  
			\draw[->, >=latex', shorten >=3pt, bend left=0,very thick] (dummy11.south) to (abs.west); 
		\end{tikzpicture}
	}
	\caption{Diagram illustrating the connections between abstract group pages and other collections of objects in the LMFDB. Nodes in solid (green) and dotted (grey) rectangles are on the production and beta version of the LMFDB, respectively.}
	\label{connectiondiagram}
\end{figure}
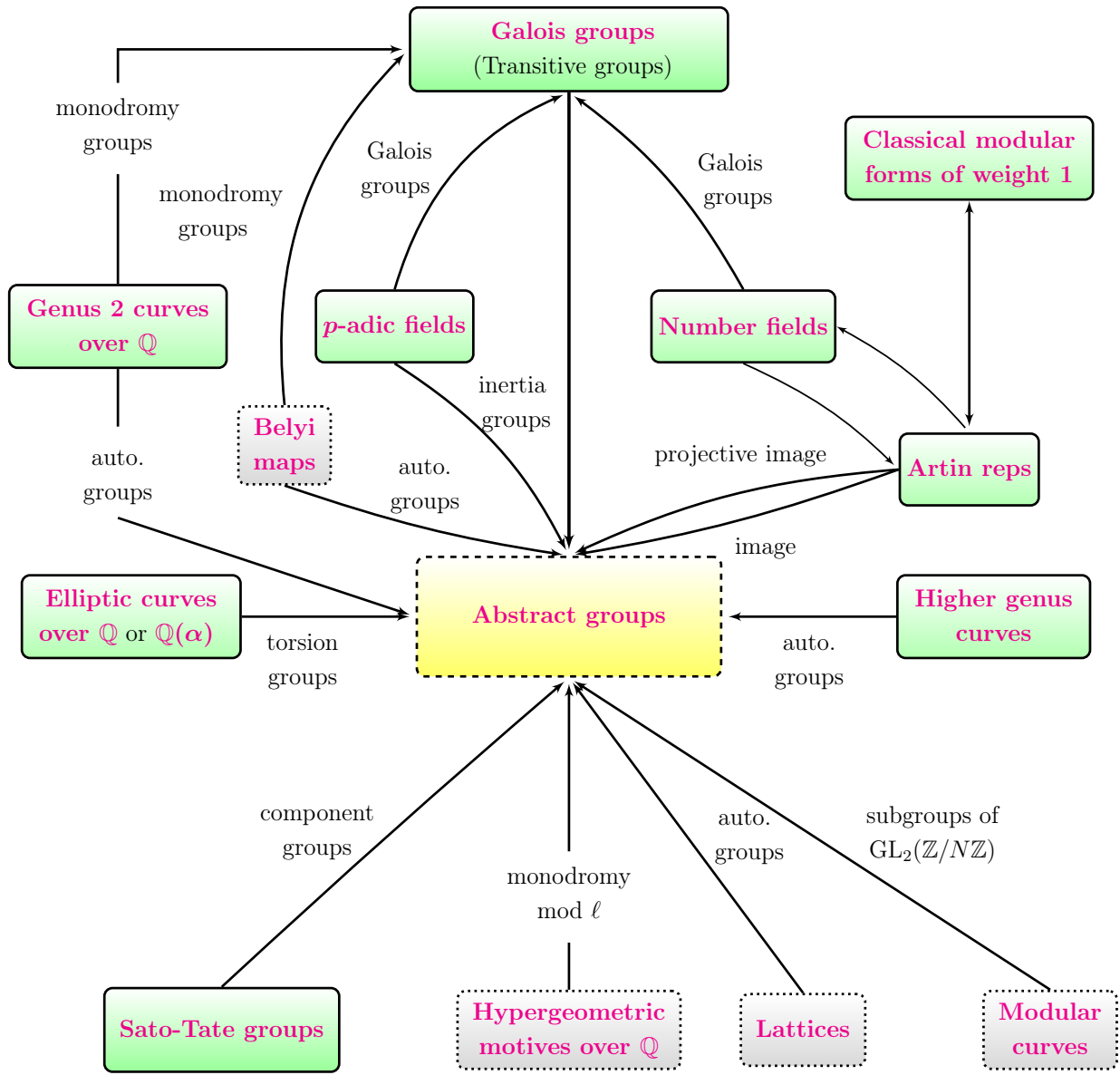

One of the integral parts of the LMFDB is demonstrating links between different mathematical objects. The \href{https://www.lmfdb.org/universe}{LMFDB universe} gives the big picture of various connections that are related to the \emph{Langlands program}. There are many related objects and links between them are displayed in the database as well. In this section we list other databases in the LMFDB which are linked with the database of Abstract Groups. We describe the various connections to other LMFDB pages with the group page in Figure~\ref{connectiondiagram}.

In Figure~\ref{connectiondiagram}, the boxes representing the LMFDB sections Belyi maps, Hypergeometric motives over $\Q$, Lattices, and Modular curves are currently only available on the beta version of the LMFDB (\url{https://beta.lmfdb.org/}).  Artin representations connect to Abstract groups via both their image and their projective images.  One can also get from Artin representations to Galois groups (i.e., transitive permutation groups) by taking a stem field for the field cut out by the image of the representation to get a number field, and then taking the Galois group of its Galois closure.  For a number field, one has the Artin representations of the Galois group of the normal closure of the field.

\begin{example}
Many objects in the LMFDB link to groups.  Having a built-in database of groups makes it easy for the user to find more information about the groups in question.  Some groups may be familiar to the user, such as those which arise as automorphism groups of genus $2$ curves; others are more complicated.  
Consider the Artin representation \url{https://www.lmfdb.org/ArtinRepresentation/4.5744.8t39.d.a} which has the smallest conductor of all $4$-dimensional irreducible Artin representations with group \url{https://www.lmfdb.org/Groups/Abstract/192.1493} \cite{jr-artin}.  This group is displayed as $C_2^3:S_4$, but since the semidirect product notation may not uniquely identify a group, the connection to a specific group page lets the user obtain lots of information about this group.  Similarly, the projective image of this group is $C_2^2:S_4$, to which the user is referred to the page \url{https://www.lmfdb.org/Groups/Abstract/96.227} for more details.
\end{example}

\subsection{Applications via Galois theory} \label{S:Galois}

Perhaps the most celebrated application of groups is their usage in studying field extensions as a part of Galois theory. Below we give examples of how our database of finite groups, together with the Galois correspondence, can be used to study objects in other collections in the LMFDB.

\begin{example}
Let $F = \Q(\alpha)$ be the number field with LMFDB label \href{https://www.lmfdb.org/NumberField/5.1.35152.1}{\texttt{5.1.35152.1}}, where $\alpha$ has minimal polynomial $f \colonequals x^5 - x^4 + 2 x^3 - 4 x^2 + x - 1$. Then $F$ has Galois closure $L = \Q(\beta)$ with LMFDB label \href{https://www.lmfdb.org/NumberField/20.0.3354518684571451850752.1}{\texttt{20.0.3354518684571451850752.1}}, where $\beta$ has minimal polynomial
\begin{align*}
x^{20} - 2 x^{19} + 10 x^{17} - 15 x^{16} &+ 40 x^{14} - 64 x^{13} + 46 x^{12} + 8 x^{11} - 32 x^{10} + 8 x^{9}\\
&+ 46 x^{8} - 64 x^{7} + 40 x^{6} - 15 x^{4} + 10 x^{3} - 2 x + 1 \, ,
\end{align*}
and $G \colonequals \Gal(L/\Q) \cong F_5$, the Frobenius group of order 20 (with LMFDB label \href{https://www.lmfdb.org/Groups/Abstract/20.3}{\texttt{20.3}}.) Examining the lattice of subgroups up to conjugacy given on the \texttt{20.3} homepage, we see that the extension $F_5$ is solvable, so the roots of $f$ can be expressed by radicals. Computing fixed fields, we obtain the subfield lattice shown in Figure \ref{fig:Galois-example},
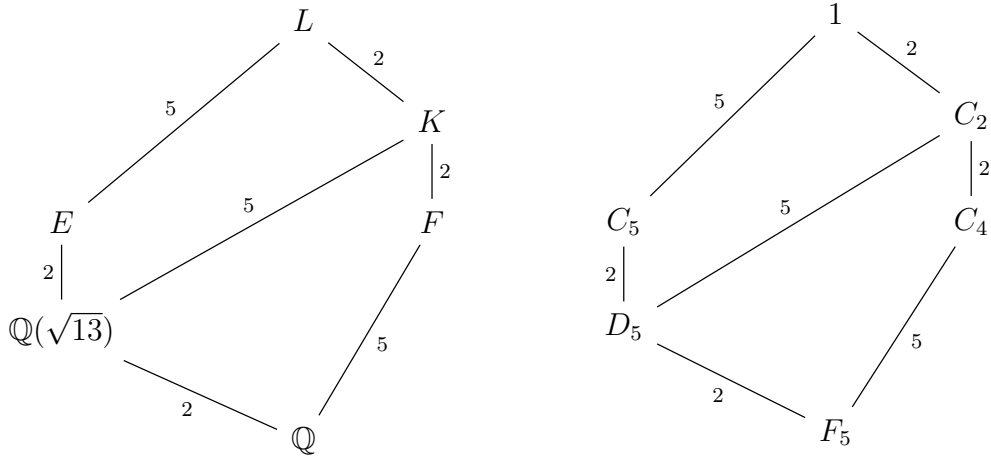
\begin{figure}[h]
    \centering
    \begin{subfigure}[h]{0.45\textwidth}
        \[\begin{tikzcd}
        	&& L \\
        	&&& {K} \\
        	{E} &&& {F} \\
        	{\Q(\sqrt{13})} \\
        	&& {\Q}
        	\arrow["5"', no head, from=1-3, to=3-1]
        	\arrow["2", no head, from=1-3, to=2-4]
        	\arrow["2", no head, from=2-4, to=3-4]
        	\arrow["5", no head, from=3-4, to=5-3]
        	\arrow["2"', no head, from=3-1, to=4-1]
        	\arrow["2"', no head, from=4-1, to=5-3]
        	\arrow["5", no head, from=4-1, to=2-4]
        \end{tikzcd}\]
    \end{subfigure}
    ~ 
    \begin{subfigure}[h]{0.45\textwidth}
        \[\begin{tikzcd}
        	&& 1 \\
        	&&& {C_2} \\
        	{C_5} &&& {C_4} \\
        	{D_5} \\
        	&& {F_5}
        	\arrow["5"', no head, from=1-3, to=3-1]
        	\arrow["2", no head, from=1-3, to=2-4]
        	\arrow["2", no head, from=2-4, to=3-4]
        	\arrow["5", no head, from=3-4, to=5-3]
        	\arrow["2"', no head, from=3-1, to=4-1]
        	\arrow["2"', no head, from=4-1, to=5-3]
        	\arrow["5", no head, from=4-1, to=2-4]
        \end{tikzcd}\]
    \end{subfigure}
    \caption{The subfield lattice of the number field with LMFDB label \href{https://www.lmfdb.org/NumberField/20.0.3354518684571451850752.1}{\texttt{20.0.3354518684571451850752.1}} and the subgroup lattice of its Galois group~$F_5$.}
    \label{fig:Galois-example}
\end{figure}
where $K$ and $E$ are the number fields given by
\begin{equation*}
    x^{10} - x^9 - 3 x^8 + 5 x^6 + x^5 - 5 x^4 + 3 x^2 - x - 1
\end{equation*}
and
\begin{equation*}
    x^4 + x^3 + 2 x^2 - 4 x + 3 \, ,
\end{equation*}
respectively. We observe that the commutator subgroup of $G$ is $C_5$, the subgroup corresponding to $E$. Thus $E$ is the largest abelian subfield of $L$, with Galois group $\Gal(E/\Q) \cong G^{\text{ab}} \cong C_4$.
\end{example}


\begin{example}
Let $K = \Q_2(\alpha)$ be the $p$-adic field with LMFDB label \href{https://www.lmfdb.org/padicField/2.4.6.7}{\texttt{2.4.6.7}}, where $\alpha$ has minimal polynomial $f \colonequals x^4 + 2x^3 + 2x^2 + 2$. Note that $K$ is not Galois over $\Q_2$ and $K$ has Galois closure $\Kgal= \Q_2(\beta)$ with LMFDB label \href{https://www.lmfdb.org/padicField/2.12.18.59}{2.12.18.59}, where $\beta$ has minimal polynomial
\begin{align*}
x^{12} - 2x^{11} + 6x^{10} + 4x^9 + 6x^8 + 12x^7 - 4x^6 - 8x^3 + 16x^2 - 8
\end{align*}
and $G \colonequals \Gal(\Kgal/\Q_2) \cong A_4$, the alternating group of order $12$ (with LMFDB label \href{https://www.lmfdb.org/Groups/Abstract/12.3}{\texttt{12.3}}).
The inertia group $I:=I(\Kgal/\Q_2)$ of $\Kgal$ is the abelian group $C_2^2$ of order $4$ (with LMFDB label \href{https://www.lmfdb.org/Groups/Abstract/4.2}{\texttt{4.2}}). The wild inertia group of $\Kgal$ (i.e., the unique $2$-Sylow subgroup of $I$) in this case is equal to $I$.

Examining the lattice of subgroups up to conjugacy given on the \href{https://www.lmfdb.org/Groups/Abstract/12.3}{\texttt{12.3}} homepage, we obtain the subfield lattice of $\Kgal/\Q_2$, shown in Figure \ref{fig:padic-example}. It is clear from the subgroup lattice of \href{https://www.lmfdb.org/Groups/Abstract/12.3}{\texttt{12.3}} that  $\Gal(\Kgal/K)\cong C_3$ and the extension $K/\Q_2$ is primitive. Also note that $A_4$ is solvable of length $2$, so we find two intermediate fields $E$ and $F$ of $\Kgal$ such that $\Gal(\Kgal/E)\cong C_2^2$ and $\Gal(\Kgal/F)\cong C_2$. 
Computing fixed fields, we obtain the subfield lattice shown in Figure \ref{fig:padic-example}, where $E$ and $F$ are the $p$-adic fields given by $ x^{3} - x + 1$ and $ x^{6} + x^{2} - 1$,
respectively. The LMFDB labels for $E$ and $F$ are \href{https://www.lmfdb.org/padicField/2.3.0.1}{\texttt{2.3.0.1}}  and \href{https://www.lmfdb.org/padicField/2.6.6.1}{\texttt{2.6.6.1}}, respectively. We observe that the commutator subgroup of $G$ is $C_2^2$, the subgroup corresponding to $E$. Thus $E$ is the largest abelian subfield of $\Kgal$, with Galois group $\Gal(E/\Q) \cong G^{\text{ab}} \cong C_3$.

This is an interesting example since $K$ is the only degree $4$ extension of $\Q_p$, for any $p$, which has Galois group $\Gal(\Kgal/\Q_2) \cong A_4$.

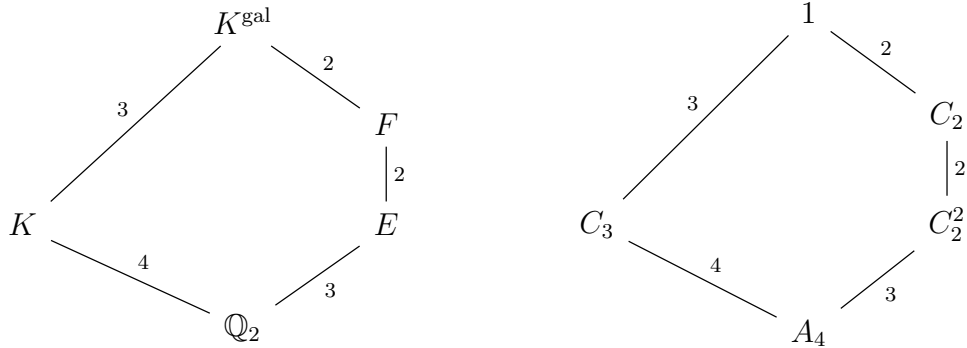
\begin{figure}[h]
	\centering
        \begin{subfigure}[H]{0.45\textwidth}
    	\[\begin{tikzcd}
        	&& \Kgal \\
        	&&& {F} \\
        	{K} &&& {E} \\
        	&& {\Q_2}
        	\arrow["3"', no head, from=1-3, to=3-1]
        	\arrow["2", no head, from=1-3, to=2-4]
        	\arrow["2", no head, from=2-4, to=3-4]
        	\arrow["3", no head, from=3-4, to=4-3]
        	\arrow["4", no head, from=3-1, to=4-3]
        \end{tikzcd}\]
    	\end{subfigure}
    	~ 
    	\begin{subfigure}[h]{0.45\textwidth}
    	\[\begin{tikzcd}
        	&& 1 \\
        	&&& {C_2} \\
        	{C_3} &&& {C_2^2} \\
        	&& {A_4}
        	\arrow["3"', no head, from=1-3, to=3-1]
        	\arrow["2", no head, from=1-3, to=2-4]
        	\arrow["2", no head, from=2-4, to=3-4]
        	\arrow["3", no head, from=3-4, to=4-3]
        	\arrow["4", no head, from=3-1, to=4-3]
        \end{tikzcd}\]
    	\end{subfigure}	
    	\caption{The subfield lattice of the Galois closure of the $p$-adic field with LMFDB label \href{https://www.lmfdb.org/padicField/2.4.6.7}{\texttt{2.4.6.7}} and the subgroup lattice of its Galois group~$A_4$.}
    	\label{fig:padic-example}
\end{figure}

\end{example}

\begin{example}
    Let $G$ be the Sato-Tate group $J(O)$ with LMFDB label \href{https://www.lmfdb.org/SatoTateGroup/1.4.F.48.48a}{\texttt{1.4.F.48.48a}}. Let $G^0$ be the connected component of the identity and $G/G^0$ be the component group of $G$. As noted on its homepage, $J(O)$ has the largest component group ($C_2 \times S_4$, LMFDB label \href{https://www.lmfdb.org/Groups/Abstract/48.48}{\texttt{48.48}}, with order 48) among Sato-Tate groups of abelian surfaces over number fields.
    
    Let $A$ be an abelian variety of dimension $g \leq 3$ defined over a number field $k$. Let $K$ be the minimal extension of $k$ over which all endomorphisms of $A$ are defined, i.e., such that $\End(A_K) = \End(A_{\overline{k}})$. By \cite[Proposition 2.17]{FKRS}, then $\ST_A/\ST_A^0 \cong \Gal(K/k)$, where $\ST_A$ is the Sato-Tate group of $A$.
    
    Let $C: y^2 = x^6 - 5x^4 + 10x^3 - 5x^2 + 2x -1$ considered over $\Q$, and let $A = \Jac(C)$. As shown in \cite{FKRS}, the variety $A$ has Sato-Tate group $J(O)$, realizing this group over $\Q$, and the endomorphisms of $A$ are defined over the number field $K = \Q(\sqrt{-2}, \sqrt{-11}, a, b)$ where
        \begin{equation*}
        a^3 - 7a + 7 = 0 \qquad \text{and} \qquad b^4 + 4b^2 + 8b + 8 = 0 \, .
    \end{equation*}
    Applying the above proposition to this example, we have
    \begin{align*}
        C_2 \times S_4 \cong G/G^0 \cong \Gal(K/k) \, .
    \end{align*}

    However, as shown in \cite[Table 4]{FKRS}, the curve $C$ can also be used to realize 24 other Sato-Tate groups by varying the base field. In other words, by taking a number field $L$ with $k \subseteq L \subseteq K$ and considering the base change $A_L$, we can obtain other Sato-Tate groups. For instance, examining \cite[Table 8]{FKRS} we see that $J(T)$ (LMFDB label \href{https://www.lmfdb.org/SatoTateGroup/1.4.F.24.13a}{\texttt{1.4.F.24.13a}}) is the unique Sato-Tate group occurring in genus $2$ with component group $C_2 \times A_4$. Computing the fixed field $L$ of $C_2 \times A_4 \leq C_2 \times S_4$, we find that $L = \Q(\sqrt{-11})$. Thus the base change $A_L$ has Sato-Tate group $J(T)$, realizing another of the 52 possible Sato-Tate groups.
    
    With further calculation, one can use other base changes of $A$ to realize other Sato-Tate groups whose component groups are subgroups of $C_2 \times S_4$ by considering the lattice of subgroups given on the homepage for $C_2 \times S_4$ and computing fixed fields.
    
\end{example}

Galois theory can also be applied in a geometric context. The equivalence of categories between function field extensions in one variable and non-singular projective curves (see \cite[\href{https://stacks.math.columbia.edu/tag/0BXX}{Tag 0BXX}]{stacks}, for instance) allows us to study non-constant morphisms of curves by examining the associated extension of function fields. We consider an example coming from the Families of higher genus curves with automorphisms collection, available at \url{https://www.lmfdb.org/HigherGenus/C/Aut/}.

\begin{example}
Consider the refined passport with label \href{https://www.lmfdb.org/HigherGenus/C/Aut/3.168-42.0.2-3-7}{\texttt{3.168-42.0.2-3-7}}. It corresponds to a unique topological equivalence class of morphisms $X \to X/\Aut(X) \cong \mathbb{P}^1$, where $X$ is the Klein quartic curve, a genus $3$ curve with the largest possible number of automorphisms for its genus. (Such curves are known as {\em Hurwitz curves}.) It has automorphism group $\Aut(X) \cong \PSL_2(\F_7)$ (which has LMFDB label \href{https://www.lmfdb.org/Groups/Abstract/168.42}{\texttt{168.42}}), and all automorphisms are defined over the field $K = \Q(\zeta_7)$. By examining the lattice of subgroups of $\PSL_2(\F_7)$ given on its homepage and applying the Galois correspondence, we can find all intermediate covers $Y$ with $X \to Y \to X/\Aut(X) \cong \PP^1$. 

Computationally, this can be accomplished using \textsf{Magma}'s \texttt{CurveQuotient} command \cite{Magma}. Calling this on each subgroup of $\PSL_2(\F_7)$ in turn, we find three intermediate covers $Y$ of genus $1$, corresponding to the subgroups of $\PSL_2(\F_7)$ isomorphic to $C_2, C_3$, and $C_4$. (All other quotients by non-trivial subgroups of $\Aut(X)$ result in genus $0$ curves.) Each of these curves can be equipped with the structure of an elliptic curve by taking as the origin of the group law the image of $(1:0:0)$ under the appropriate quotient map. With some further computation (see \cite[equation 2.10]{Elkies}), one can show that each of these elliptic curves is isomorphic to the curve $E: y^2 = 4 x^3 + 21 x^2 + 28 x$ over $K$. From this we observe that $\Jac(X)$ decomposes as $E^3$ up to isogeny.
\end{example}

\subsection{Modular curves and subgroups of \texorpdfstring{$\GL_2(\Z/N\Z)$}{GL(2,Z/NZ)}}
The LMFDB currently has a preliminary database of modular curves, available at \url{https://www.lmfdb.org/ModularCurve/Q/}. We conclude this section by describing how properties of modular curves can be deduced from their corresponding subgroups of $\GL_2(\Z/N\Z)$. For a more comprehensive exposition of modular curves and subgroups of $\GL_2(\Zhat)$, see \cite{RSZB} or \cite{Zywina}.

Let $E$ be an elliptic curve over $\Q$ and $E[N]$ be its $N$-torsion subgroup for each $N \in \Z_{\geq 1}$. The absolute Galois group $G_\Q \colonequals \Gal(\overline{\Q}/\Q)$ acts on $E[N]$, and since $E[N] \cong (\Z/N\Z)^2$ as abelian groups, we obtain a representation
\begin{equation*}
    \rho_{E,N} : G_\Q \to \Aut(E[N]) \cong \GL_2(\Z/N\Z) \, .
\end{equation*}
By choosing compatible bases, we can take the inverse limit and package these together as a single representation
\begin{equation*}
    \rho_E: G_\Q \to \varprojlim_N \GL_2(\Z/N\Z) = \GL_2(\Zhat) \, .
\end{equation*}
If $E$ does not have complex multiplication, then Serre's Open Image Theorem \cite{Serre} implies that the image of $\rho_E$ is an open subgroup of $\GL_2(\Zhat)$, hence has finite index. Given an open subgroup $H$ of $\GL_2(\Zhat)$, one can define the modular curve $X_H$ whose $K$-points parametrize elliptic curves $E/K$ such that $\img(\rho_E) \subseteq H$ (up to conjugation). (For a precise definition of $X_H$, see \cite[\S2.3]{RSZB} or \cite[\S3]{Zywina}.)

For each $N \in \Z_{\geq 1}$, let $\pi_N : \GL_2(\Zhat) \to \GL_2(\Z/N\Z)$ be the projection map. Every open subgroup $H \leq \GL_2(\Zhat)$ contains $\ker(\pi_N)$ for some $N$, and the smallest such $N \in \Z_{\geq 1}$ is called the \defn{level} of $H$. If $H$ contains $\ker(\pi_N)$, then $H$ can be recovered from $\pi_N(H)$ as $H = \pi_N^{-1}(\pi_N(H))$. Thus to store $H$ on a computer, we can simply store generators for $\pi_N(H)$ where $N$ is the level of $H$.

 Although this is a database of modular curves, much of the geometric data is computed group theoretically.
 Let $H \leq \GL_2(\Zhat)$ be an open subgroup of level $N$. Letting $\Gamma_H \colonequals \pm \pi_N(H) \cap \SL_2(\Z/N\Z)$, then one can determine the number of elliptic points and cusps of $X_H$ by studying the action of the matrices
\begin{equation*}
    \begin{pmatrix} 0 & 1\\ -1 & 0 \end{pmatrix}, \quad
    \begin{pmatrix} 0 & 1\\ -1 & -1 \end{pmatrix}, \quad
    \begin{pmatrix} 1 & 1\\ 0 & 1 \end{pmatrix}
\end{equation*}
on the coset space $\Gamma_H \backslash \SL_2(\Z/N\Z)$. This data can in turn be used to compute the genus of $X_H$.
    
\begin{example}
    Consider the modular curve $X_0(6)$ of level $6$ with LMFDB label \href{https://www.lmfdb.org/ModularCurve/Q/6.12.0.a.1/}{\texttt{6.12.0.a.1}}. As a containment $H \subseteq H'$ of open subgroups of $\GL_2(\Zhat)$ induces a morphism $X_H \to X_{H'}$ of modular curves, we have morphisms $X_0(6) \to X_0(2)$ and $X_0(6) \to X_0(3)$. On the homepage for $X_0(6)$, it is further claimed that $X_0(6)$ is the fiber product of $X_0(2)$ and $X_0(3)$ over $X(1)$---we explain how this can be determined group theoretically.
    
    Let $H_2, H_3, H_6$ be the subgroups of $\GL_2(\Z/2\Z)$, $\GL_2(\Z/3\Z)$, $\GL_2(\Z/6\Z)$ corresponding to $X_0(2)$ (with label \href{https://www.lmfdb.org/ModularCurve/Q/2.3.0.a.1/}{\texttt{2.3.0.a.1}}), $X_0(3)$ (with label \href{https://www.lmfdb.org/ModularCurve/Q/3.4.0.a.1/}{\texttt{3.4.0.a.1}}), and $X_0(6)$, respectively. Note that $\GL_2(\Z/6\Z)$ has label \href{https://www.lmfdb.org/Groups/Abstract/288.851}{\texttt{288.851}} and the subgroup $H_6 \cong C_2^2 \times S_3$ has  subgroup label \href{https://www.lmfdb.org/Groups/Abstract/sub/288.851.12.c1.a1}{\texttt{288.851.12.c1.a1}}.
    Taking the inverse image of $H_2$ and $H_3$ under the projection maps $\pi_{6,2}: \GL_2(\Z/6\Z) \to \GL_2(\Z/2\Z)$, $\pi_{6,3}: \GL_2(\Z/6\Z) \to \GL_2(\Z/3\Z)$, we obtain subgroups $\wt{H_2} \colonequals \pi_{6,2}^{-1}(H_2)$ with subgroup label
    \href{https://www.lmfdb.org/Groups/Abstract/sub/288.851.3.a1.a1}{\texttt{288.851.3.a1.a1}}
    and $\wt{H_3} \colonequals \pi_{6,3}^{-1}(H_3)$ with subgroup label \href{https://www.lmfdb.org/Groups/Abstract/sub/288.851.4.b1.a1}{\texttt{288.851.4.b1.a1}}.
    
    Either by examining the (very large) \href{https://www.lmfdb.org/Groups/Abstract/diagram/288.851}{subgroup lattice diagram} of $\GL_2(\Z/6\Z)$ or by computing directly, we find that $\wt{H_2} \cap \wt{H_3} = H_6$ and $\langle \wt{H_2}, \wt{H_3} \rangle = \GL_2(\Z/6\Z)$. This shows that $H_6$ is the fiber product of $\wt{H_2}$ and $\wt{H_3}$ over $\GL_2(\Z/6\Z)$, and hence $X_0(6)$ is the fiber product of $X_0(2)$ and $X_0(3)$ over $X(1)$, as depicted in Figure \ref{fig:fiber-product-example}.
    
    \begin{figure}[h]
    \centering
    \begin{subfigure}[h]{0.45\textwidth}
        \[\begin{tikzcd}
        	{X_0(6)} \arrow[dr, phantom, "\text{\textopencorner}", pos=0.2] & {X_0(2)} \\
        	{X_0(3)} & {X(1)}
        	\arrow[from=1-2, to=2-2]
        	\arrow[from=2-1, to=2-2]
        	\arrow[from=1-1, to=2-1]
        	\arrow[from=1-1, to=1-2]
        \end{tikzcd}\]
    \end{subfigure}
    ~ 
    \begin{subfigure}[h]{0.45\textwidth}
        \[\begin{tikzcd}
        	{H_6} \arrow[dr, phantom, "\text{\textopencorner}", pos=0.2] & {\wt{H_2}} \\
        	{\wt{H_3}} & {\GL_2(\Z/6\Z)}
        	\arrow[from=1-2, to=2-2]
        	\arrow[from=2-1, to=2-2]
        	\arrow[from=1-1, to=2-1]
        	\arrow[from=1-1, to=1-2]
        \end{tikzcd}\]
    \end{subfigure}
    \caption{The fiber product diagram for the modular curve $X_0(6)$, and the fiber product diagram of the corresponding subgroups of $\GL_2(\Z/6\Z)$.}
    \label{fig:fiber-product-example}
\end{figure}
    
\end{example}




\bibliographystyle{alpha} 
\bibliography{references}

\end{document}